\documentclass{amsart}
\usepackage{mathrsfs,url,amssymb}
\newtheorem{thm}{Theorem}[section]
\newtheorem{cor}[thm]{Corollary}
\newtheorem{lem}[thm]{Lemma}
\newtheorem{prop}[thm]{Proposition}
\theoremstyle{definition}
\newtheorem{defn}[thm]{Definition}
\newtheorem{fact}[thm]{Fact}
\newtheorem{rem}[thm]{Remark}
\newtheorem{que}[thm]{Question}

\numberwithin{equation}{section}

\begin{document}
\subjclass[2000]{Primary 20E15; Secondary 06E10, 20A15}


\title{Strongly bounded groups and infinite powers of finite groups}
\author{Yves de Cornulier}%
\date{\today}
\begin{abstract}
We define a group as strongly bounded if every isometric action on
a metric space has bounded orbits. This latter property is
equivalent to the so-called uncountable strong cofinality,
recently introduced by Bergman.

Our main result is that $G^I$ is strongly bounded when $G$ is a
finite, perfect group and $I$ is any set. This strengthens a
result of Koppelberg and Tits. We also prove that
$\omega_1$-existentially closed groups are strongly bounded.
\end{abstract}
\maketitle
\section{Introduction}
Let us say that a group is strongly bounded if every isometric
action on a metric space has bounded orbits.

We observe that the class of discrete, strongly bounded groups
coincides with a class of groups which has recently emerged since
a preprint of Bergman \cite{Bergman}, sometimes referred to as
``groups with uncountable strong cofinality", or ``groups with
Bergman's Property". This class contains no countably infinite
group, but contains symmetric groups over infinite sets
\cite{Bergman}, various automorphism groups of infinite structures
such as 2-transitive chains \cite{DH}, full groups of certain
equivalence relations \cite{Miller}, oligomorphic permutation
groups with ample generics \cite{KR}; see \cite{Bergman} for more
references.

In Section \ref{s:exist}, we prove that $\omega_1$-existentially
closed groups are strongly bounded. This strengthens a result of
Sabbagh \cite{Sabbagh}, who proved that they have cofinality
$\neq\nolinebreak\omega$.

In Section \ref{s:power}, we prove that if $G$ is any finite
perfect group, and $I$ is any set, then $G^I$, endowed with the
discrete topology, is strongly bounded. This strengthens a result
of Koppelberg and Tits \cite{KT}, who proved that this group has
Serre's Property (FA). This group has finite exponent and is
locally finite, hence amenable. In contrast, all previously known
infinite strongly bounded groups contain a non-abelian free
subgroup.

\section{Strongly bounded groups}

\begin{defn}
We say that a group $G$ is {\em strongly bounded} if every
isometric action of $G$ on a metric space has bounded orbits.
\end{defn}

\begin{rem}
Let $G$ be a strongly bounded group. Then every isometric action
of $G$ on a nonempty complete $\textnormal{CAT}(0)$ space has a
fixed point; in particular, $G$ has Property (FH) and Property
(FA), which mean, respectively, that every isometric action of $G$
on a Hilbert space (resp. simplicial tree) has a fixed point. This
follows from the Bruhat-Tits fixed point lemma, which states that
every action of a group on a complete $\textnormal{CAT}(0)$ space
which has a bounded orbit has a fixed point \cite[Chap. II,
Corollary 2.8(1)]{BH}.

It was asked in \cite{W} whether the equivalence between Kazhdan's
Property (T) and Property (FH), due to Delorme and Guichardet (see
\cite[Chap. 2]{BHV}) holds for more general classes of groups than
locally compact $\sigma $-compact groups; in particular, whether
it holds for general locally compact groups.

The answer is negative, even if we restrict to discrete groups:
this follows from the existence of uncountable strongly bounded
groups, combined with the fact that Kazhdan's Property (T) implies
finite generation \cite[Chap. 1]{BHV}.\label{rem:sb_FH}\end{rem}

\begin{defn}
We say that a group $G$ is {\em Cayley bounded} if, for every
generating subset $U\subseteq G$, there exists some $n$ (depending
on $U$) such that every element of $G$ is a product of $n$
elements of $U\cup U^{-1}\cup\{1\}$. This means every Cayley graph
of $G$ is bounded.

A group $G$ is said to have cofinality $\omega$ if it can be
expressed as the union of an increasing sequence of proper
subgroups; otherwise it is said to have cofinality $\neq\omega$.
\label{def_strongly_bounded}
\end{defn}

The combination of these two properties, sometimes referred as
``uncountable strong cofinality"\footnote{In the literature, it is
sometimes referred as ``Bergman's Property"; Bergman's Property
also sometimes refers to Cayley boundedness without cofinality
assumption.}, has been introduced and is extensively studied in
Bergman's preprint \cite{Bergman}; see also \cite{DG}. Note that
an uncountable group with cofinality $\neq\omega$ is not
necessarily Cayley bounded: the free product of two uncountable
groups of cofinality $\neq\omega$, or the direct product of an
uncountable group of cofinality $\neq\omega$ with $\mathbf{Z}$,
are obvious counterexamples. On the other hand, a Cayley bounded
group with cofinality $\omega$ is announced in \cite{Khelif}.

\medskip

The following result can be compared to \cite[Lemma 10]{Bergman}:

\begin{prop}
A group $G$ is strongly bounded if and only if it is Cayley
bounded and has cofinality~$\neq\omega$.
\label{strongbd}
\end{prop}
\noindent \textbf{Proof}:~Suppose that $G$ is not Cayley bounded.
Let $U$ be a generating subset such that $G$ the corresponding
Cayley graph is not bounded. Since $G$ acts transitively on it, it
has an unbounded orbit.

Suppose that $G$ has cofinality $\omega$. Then $G$ acts on a tree
with unbounded orbits \cite[Chap I, \S 6.1]{Serre}.

\medskip

Conversely, suppose that $G$ has has cofinality $\neq\omega$ and
is Cayley bounded. Let $G$ act isometrically on a metric space.
Let $x\in X$, let $K_n=\{g\in G\,|\;d(x,gx)<n\}$, and let $H_n$ be
the subgroup generated by $K_n$. Then $G=\bigcup K_n=\bigcup H_n$.
Since $G$ has cofinality $\neq\omega$, $H_n=G$ for some $n$, so
that $K_n$ generates $G$. Since $G$ is Cayley bounded, and since
$K_n$ is symmetric, $G\subseteq (K_n)^m$ for some $m$. This easily
implies that $G\subseteq K_{nm}$, so that the orbit of $x$ is
bounded.~$\blacksquare$

\begin{rem}
It follows that a countably infinite group $\Gamma$ is not
strongly bounded: indeed, either $\Gamma$ has a finite generating
subset, so that the corresponding Cayley graph is unbounded, or
else $\Gamma$ is not finitely generated, so is an increasing union
of a sequence of finitely generated subgroups, so has
cofinality~$\omega$.
\end{rem}

\medskip

\begin{defn}If $G$ is a group, and $X\subseteq G$, define
$$\mathcal{G}(X)=X\cup\{1\}\cup\{x^{-1},\;x\in
X\}\cup\{xy\,|\;x,y\in X\}.$$\label{G(X)}\end{defn}

The following proposition is immediate and is essentially
contained in Lemma 10 of \cite{Bergman}.

\begin{prop}
The group $G$ is strongly bounded if and only if, for every
increasing sequence $(X_n)$ of subsets such that $\bigcup_n X_n=G$
and $\mathcal{G}(X_{n})\subseteq X_{n+1}$ for all $n$, one has
$X_n=G$ for some~$n$.~$\blacksquare$\label{strong_bound}
\end{prop}

\begin{rem}
The first Cayley bounded groups with uncountable cofinality were
constructed by Shelah \cite[Theorem 2.1]{Shelah}. They seem to be
the only known to have a uniform bound on the diameter of Cayley
graphs. They are torsion-free. These groups are highly
non-explicit and their construction, which involves small
cancellation theory, rests on the Axiom of Choice.

The first explicit examples, namely, symmetric groups over
infinite sets, are due to Bergman \cite{Bergman}. The first
explicit torsion-free examples, namely, automorphism groups of
double transitive chains, are due to Droste and Holland \cite{DH}.
\end{rem}

\begin{rem}
It is easy to observe that groups with cofinality $\neq\omega$
also have a geometric characterization; namely, a group $G$ has
cofinality $\neq\omega$ if and only if every isometric action of
$G$ on an \textit{ultrametric} metric space has bounded orbits.
\end{rem}

\begin{rem}
In \cite[\S 2.6]{BHV}, it is proved that a infinite solvable group
never has Property~(FH) (defined in Remark \ref{rem:sb_FH}). In
particular, an infinite solvable group is never strongly bounded.
This latter result is improved by Khelif \cite{Khelif} who proves
that an infinite solvable group is never Cayley bounded. On the
other hand, it is not known whether there exist uncountable
solvable groups with cofinality~$\neq\omega$.\label{rem:solv_SB}
\end{rem}

\section{$\omega_1$-existentially closed groups}\label{s:exist}

Recall that a group $G$ is $\omega_1$-existentially closed if
every countable set of equations and inequations with coefficients
in $G$ which has a solution in a group containing $G$, has a
solution in $G$. Sabbagh \cite{Sabbagh} proved that every
$\omega_1$-existentially closed group has cofinality $\neq\omega$.
We give a stronger result, which has been independently noticed by
Khelif \cite{Khelif}:

\begin{thm}
Every $\omega_1$-existentially closed group $G$ is strongly
bounded.
\end{thm}
\noindent \textbf{Proof}:~Let $G$ act isometrically on a nonempty
metric space $X$. Fix $x\in X$, and define $\ell(g)=d(gx,x)$ for
all $g\in G$. Then $\ell$ is a length function, i.e. satisfies
$\ell(1)=0$ and $\ell(gh)\le \ell(g)+\ell(h)$ for all $g,h\in G$.
Suppose by contradiction that $\ell$ is not bounded. For every
$n$, fix $c_n\in G$ such that $\ell(c_n)\ge n^2$. Let $C$ be the
group generated by all $c_n$. By the proof of the HNN embedding
Theorem \cite[Theorem 3.1]{LS}, $C$ embeds naturally in the group
$$\Gamma=\langle C,a,b,t\;;\;
c_n=t^{-1}b^{-n}ab^nta^{-n}b^{-1}a^n\;(n\in\mathbf{N})\rangle,$$
which is generated by $a,b,t$. Since $G$ is
$\omega_1$-existentially closed, there exist
$\bar{a},\bar{b},\bar{t}$ in $G$ such that the group generated by
$C$, $\bar{a}$, $\bar{b}$, and $\bar{t}$ is naturally isomorphic
to $\Gamma$. Set
$M=\max(\ell(\bar{a}),\ell(\bar{b}),\ell(\bar{t}))$. Then, since
$\ell$ is a length function and $c_n$ can be expressed by a word
of length $4n+4$ in $a,b,c$, we get $\ell(c_n)\le M(4n+4)$ for all
$n$, contradicting $\ell(c_n)\ge n^2$.~$\blacksquare$

\medskip

It is known \cite{Scott} that every group embeds in a
$\omega_1$-existentially closed group. Thus, we obtain:

\begin{cor}
Every group embeds in a strongly bounded group.~$\blacksquare$
\end{cor}

Note that this was already a consequence of the strong boundedness
of symmetric groups \cite{Bergman}, but provides a better
cardinality: if $|G|=\kappa$, we obtain a group of cardinality
$\kappa^{\aleph_0}$ rather that $2^{\kappa}$.

\section{Powers of finite groups}\label{s:power}

\begin{thm}
Let $G$ be a finite perfect group, and $I$ a set. Then the
(unrestricted) product $G^I$ is strongly
bounded.\label{G^I_strongly_bounded}
\end{thm}

\begin{rem}
Conversely, if $I$ is infinite and $G$ is not perfect, then $G^I$
maps onto the direct sum $\mathbf{Z}/p\mathbf{Z}^{(\mathbf{N})}$
for some prime $p$, so has cofinality $\omega$ and is not Cayley
bounded, as we see by taking as generating subset the canonical
basis of $\mathbf{Z}/p\mathbf{Z}^{(\mathbf{N})}$.
\end{rem}

\begin{rem}
By Theorem \ref{G^I_strongly_bounded}, every Cayley graph of $G^I$
is bounded. If $I$ is infinite and $G\neq 1$, one can ask whether
we can choose a bound which does not depend on the choice of the
Cayley graph. The answer is negative: indeed, for all
$n\in\mathbf{N}$, observe that the Cayley graph of $G^n$ has
diameter exactly $n$ if we choose the union of all factors as
generating set. By taking a morphism of $G^I$ onto $G^n$ and
taking the preimage of this generating set, we obtain a Cayley
graph for $G^I$ whose diameter is exactly $n$.
\end{rem}

Our remaining task is to prove Theorem \ref{G^I_strongly_bounded}.
The proof is an adequate modification of the original proof of the
(weaker) result of Koppelberg and Tits \cite{KT}, which states
that $G^I$ has cofinality $\neq\omega$.

\medskip

If $A$ is a ring with unity, and $X\subseteq A$, define
$$\mathcal{R}(X)=X\cup\{-1,0,1\}\cup\{x+y\,|\;x,y\in X\}\cup\{xy\,|\;x,y\in X\}.$$
It is clear that $\bigcup_{n\in\mathbf{N}}\mathcal{R}^n(X)$ is the
subring generated by $X$.

Recall that a Boolean algebra is an associative ring with unity
which satisfies $x^2=x$ for all~$x$. Such a ring has
characteristic 2 (since $2=2^2-2$) and is commutative (since
$xy-yx=(x+y)^2-(x+y)$). The ring $\mathbf{Z}/2\mathbf{Z}$ is a
Boolean algebra, and so are all its powers
$\mathbf{Z}/2\mathbf{Z}^E=\mathscr{P}(E)$, for any set $E$.

\begin{prop}
Let $E$ be a set, and $(\mathscr{X}_i)_{i\in\mathbf{N}}$ an
increasing sequence of subsets of $\mathscr{P}(E)$. Suppose that
$\mathcal{R}(\mathscr{X}_i)\subseteq \mathscr{X}_{i+1}$ for all
$i$. Suppose that
$\mathscr{P}(E)=\bigcup_{i\in\mathbf{N}}\mathscr{X}_i$. Then
$\mathscr{P}(E)=\mathscr{X}_i$ for some $i$.\label{cof_P(I)}
\end{prop}

\begin{rem}
1) We could have defined, in analogy of Definition
\ref{def_strongly_bounded}, the notion of strongly bounded ring
(although the terminology ``uncountable strong cofinality" seems
more appropriate in this context). Then Proposition \ref{cof_P(I)}
can be stated as: if $E$ is infinite, the ring
$\mathscr{P}(E)=\mathbf{Z}/2\mathbf{Z}^E$ is strongly bounded. If
$E$ is infinite, note that, as an {\em additive group}, it maps
onto $\mathbf{Z}/2\mathbf{Z}^{(\mathbf{N})}$, so has cofinality
$\omega$ and is not Cayley bounded.
\end{rem}

\noindent \textbf{Proof} of Proposition \ref{cof_P(I)}. Suppose
the contrary. If $X\subseteq E$, denote by $\mathscr{P}(X)$ the
power set of $X$, and view it as a subset of $\mathscr{P}(E)$.
Define $\mathscr{L}=\{X\in\mathscr{P}(E)\,|\;\forall
i,\;\mathscr{P}(X)\nsubseteq \mathscr{X}_i\}$. The assumption is
then: $E\in\mathscr{L}$.

Observation: if $X\in \mathscr{L}$ and $X'\subseteq X$, then
either $X'$ or $X-X'$ belongs to $\mathscr{L}$. Indeed, otherwise,
some $\mathscr{X}_i$ would contain $\mathscr{P}(X')$ and
$\mathscr{P}(X-X')$, and then $\mathscr{X}_{i+1}$ would contain
$\mathscr{P}(X)$.

We define inductively a decreasing sequence of subsets
$B_i\in\mathscr{L}$, and a non-decreasing sequence of integers
$(n_i)$ by:

$$B_0=E;$$
$$n_i=\inf\{t\,|\; B_i\in\mathscr{X}_t\};$$
$$B'_{i+1}\subset B_i\quad\text{and}\quad
B'_{i+1}\notin\mathscr{X}_{n_i+1};$$
$$B_{i+1}=\left\{%
\begin{array}{ll}
    B'_{i+1}, & \hbox{if } B'_{i+1}\in\mathscr{L}, \\
    B_i-B'_{i+1}, & \hbox{otherwise.} \\
\end{array}%
\right.$$

Define also $C_i=B_i-B_{i+1}$. The sets $C_i$ are pairwise
disjoint.

\begin{fact}For all $i$, $B_{i+1}\notin\mathscr{X}_{n_i}$ and $C_{i}\notin\mathscr{X}_{n_i}$.\end{fact}
\noindent \textit{Proof}:~Observe that
$\{B_{i+1},C_i\}=\{B'_{i+1},B_i-B'_{i+1}\}$. We already know
$B'_{i+1}\notin\mathscr{X}_{n_i+1}$, so it suffices to check
$B_i-B'_{i+1}\notin\mathscr{X}_{n_i}$. Otherwise,
$B'_{i+1}=B_i-(B_i-B'_{i+1})\in\mathcal{R}(\{B_i,B_i-B'_{i+1}\})
\subseteq\mathcal{R}(\mathscr{X}_{n_i})\subseteq\mathscr{X}_{n_i+1}$;
this is a contradiction.~$\Box$

\medskip

This fact implies that the sequence $(n_i)$ is strictly
increasing. We now use a diagonal argument. Let
$(N_j)_{j\in\mathbf{N}}$ be a partition of $\mathbf{N}$ into
infinite subsets. Set $D_j=\bigsqcup_{i\in N_j}C_i$ and
$m_j=\inf\{t\,|\;D_j\in\mathscr{X}_t\}$, and let $l_j$ be an
element of $N_j$ such that $l_j>\max(m_j,j)$.

Set $X=\bigsqcup_j C_{l_j}$. For all $j$, $D_j\cap
X=C_{l_j}\notin\mathscr{X}_{l_j}$. On the other hand,
$D_j\in\mathscr{X}_{m_j}\subseteq\mathscr{X}_{l_j-1}$ since
$l_j\ge m_j+1$. This implies $X\notin \mathscr{X}_{l_j-1}\supseteq
\mathscr{X}_j$ for all $j$, contradicting
$\mathscr{P}(E)=\bigcup_{i\in\mathbf{N}}\mathscr{X}_i$.~$\blacksquare$

\medskip

The following corollary, of independent interest, was suggested to
me by Romain Tessera.

\begin{cor}
Let $A$ be a finite ring with unity (but not necessarily
associative or commutative). Let $E$ be a set, and
$(\mathscr{X}_i)_{i\in\mathbf{N}}$ an increasing sequence of
subsets of $A^E$. Suppose that
$\mathcal{R}(\mathscr{X}_i)\subseteq \mathscr{X}_{i+1}$ for all
$i$. Suppose that $A^E=\bigcup_{i\in\mathbf{N}}\mathscr{X}_i$.
Then $A^E=\mathscr{X}_i$ for some $i$.
\end{cor}
\noindent \textbf{Proof}:~By reindexing, we can suppose that
$\mathscr{X}_0$ contains the constants. Write
$\mathscr{Y}_i=\{J\subseteq E\,|\; 1_J\in\mathscr{X}_{3i}\}$. If
$J,K\in\mathscr{Y}_i$, $1_{J\cap
K}=1_J1_K\in\mathscr{X}_{3i+1}\subseteq\mathscr{X}_{3i+3}$, so
that $J\cap K\in\mathscr{Y}_{i+1}$, and $1_{J\vartriangle
K}=1_J+1_K-2.1_J1_K\in\mathscr{X}_{3i+3}$, so that $J\vartriangle
K\in\mathscr{Y}_{i+1}$. By Proposition \ref{cof_P(I)},
$\mathscr{Y}_m=\mathscr{P}(E)$ for some $m$. It is then clear that
$A^E=\mathscr{X}_n$ for some $n$ (say,
$n=3m+1+\lceil\log_2|A|\rceil$).~$\blacksquare$

\bigskip

If $A$ is a Boolean algebra, and $X\subseteq A$, we define
$$\mathcal{D}(X)=X\cup\{0,1\}\cup\{x+y\,|\;x,y\in X\text{ such that }xy=0\}\cup\{xy\,|\;x,y\in
X\}.$$
$$\mathcal{I}_k(X)=\{x_1x_2\dots x_k\,|\; x_1,\dots,x_k\in X\}.$$
$$\mathcal{V}_k(X)=\{x_1+x_2+\dots x_k\,|\;x_1,\dots,x_k\in X
\text{ such that }x_ix_j=0\;\;\forall i\neq j\}.$$

The following lemma contains some immediate facts which will be
useful in the proof of the main result.

\begin{lem}
Let $A$ be a Boolean algebra, and $X\subseteq A$ a symmetric
subset (i.e. closed under $x\mapsto 1-x$) such that $0\in X$.
Then, for all $n\ge 0$,

1) $\mathcal{R}^n(X)\subseteq \mathcal{D}^{2n}(X)$, and

2) $\mathcal{D}^{n}(X)\subseteq
\mathcal{V}_{2^{2^n}}(\mathcal{I}_{2^n}(X))$.\label{RnD2n}
\end{lem}
\noindent \textbf{Proof}:~1) It suffices to prove
$\mathcal{R}(X)\subseteq\mathcal{D}^2(X)$. Then the statement of
the lemma follows by induction. Let $u\in\mathcal{R}(X)$. If
$u\notin\mathcal{D}(X)$, then $u=x+y$ for some $x,y\in X$. Then
$u=(1-x)y+(1-y)x\in\mathcal{D}^2(X)$.

2) Is an immediate induction.~$\blacksquare$

\begin{defn}[\cite{KT}]
Take $n\in\mathbf{N}$, and let $G$ be a group. Consider the set of
functions $G^n\to G$; this is a group under pointwise
multiplication. The elements $m(g_1,\dots,g_n)$ in the subgroup
generated by the constants and the canonical projections are
called {\em monomials}. Such a monomial is {\em homogeneous} if
$m(g_1,\dots,g_n)=1$ whenever at least one $g_i$ is equal to $1$.
\end{defn}

\begin{lem}[\cite{KT}]
Let $G$ be a finite group which is not nilpotent. Then there exist
$a\in G$, $b\in G-\{1\}$, and a homogeneous monomial $f:G^2\to G$,
such that $f(a,b)=b$.\label{f(a,b)}
\end{lem}


The proof can be found in \cite{KT}, but, for the convenience of
the reader, we have included the proof from \cite{KT} in the
(provisional) Appendix below.


\begin{rem}
If $G$ is a group, and $f(x_1,\dots,x_n)$ is a homogeneous
monomial with $n\ge 2$, then $m(g_1,\dots,g_n)=1$ whenever at
least one $g_i$ is central: indeed, we can then write, for all
$x_1,\dots,x_n$ with $x_i$ central,
$m(x_1,\dots,x_i,\dots,x_n)=m'(x_1,\dots,\widehat{x_i},\dots,x_n)x_i^k$.
By homogeneity in $x_i$,
$m'(x_1,\dots,\widehat{x_i},\dots,x_n)=1$, and we conclude by
homogeneity in $x_j$ for any $j\neq i$.

Accordingly, if $(C_\alpha)$ denotes the (transfinite) ascending
central series of $G$, an immediate induction on $\alpha$ shows
that if $f(a,b)=b$ for some homogeneous monomial $f$, $a\in G$ and
$b\in C_\alpha$, then $b=1$. In particular, if $G$ is nilpotent
(or even residually nilpotent), then the conclusion of Lemma
\ref{f(a,b)} is always false.
\end{rem}

\medskip

\begin{lem}
Let $G$ be a finite group, $I$ a set, and $H=G^I$. Suppose that
$f(a,b)=b$ for some $a,b\in G$, and some homogeneous monomial $f$,
and let $N$ be the normal subgroup of $G$ generated by $b$. Let
$(X_m)$ be an increasing sequence of subsets of $H$ such that
$\mathcal{G}(X_m)\subseteq X_{m+1}$ (see Definition \ref{G(X)}),
and $\bigcup X_m=H$. Then $N^I\subseteq X_m$ for $m$ big
enough.\label{Main_Lemma}
\end{lem}

\noindent \textbf{Proof}:~Suppose the contrary. If $x\in G$ and
$J\subseteq I$, denote by $x_J$ the element of $G^I$ defined by
$x_J(i)=x$ if $i\in J$ and $x_J(i)=1$ if $i\notin J$.

Denote by $\bar{f}=f^I$ the corresponding homogeneous monomial:
$H^2\to H$. Upon extracting, we can suppose that all $c_I$, $c\in
G$, are contained in $X_0$. In particular, the ``constants" which
appear in $\bar{f}$ are all contained in $X_0$.

Hence we have, for all $m$, $\bar{f}(X_m,X_m)\subseteq X_{m+d}$,
where $d$ depends only on the length of $f$. For $J,K\subseteq I$,
we have the following relations:

\begin{equation}
a_I.a_J^{-1}=a_{I-J},\label{eq1}
\end{equation}
\begin{equation}
\bar{f}(a_J,b_K)=b_{J\cap K},\label{eq2}
\end{equation}
\begin{equation}
\bar{f}(a_J,b_I)=b_J,\label{eq3}
\end{equation}
\begin{equation}
 \text{If }J\cap
K=\emptyset,\quad b_J.\,b_K=b_{J\sqcup K}.\label{eq4}
\end{equation}

For all $m$, write $\mathscr{W}_m=\{J\in\mathscr{P}(I)\,|\;a_J\in
X_m\}$, and let $\mathscr{A}_m$ be the Boolean algebra generated
by $\mathscr{W}_m$. Then $\bigcup_m\mathscr{A}_m=\mathscr{P}(I)$.
By Proposition \ref{cof_P(I)}, there exists some $M$ such that
$\mathscr{A}_M=\mathscr{P}(I)$. Set
$\mathscr{X}_n=\mathcal{R}^n(\mathscr{W}_M)$. Then, since
$\mathscr{A}_M=\mathscr{P}(I)$,
$\bigcup_n\mathscr{X}_n=\mathscr{P}(I)$. Again by Proposition
\ref{cof_P(I)}, there exists some $N$ such that
$\mathscr{X}_N=\mathscr{P}(I)$. So, by 1) of Lemma \ref{RnD2n}, we
get

\begin{equation}\mathcal{D}^{2N}(\mathscr{W}_M)=\mathscr{P}(I).\label{Eq5}\end{equation}

Define, for all $m$,
$\mathscr{Y}_m=\{J\in\mathscr{P}(I)\,|\;b_J\in X_m\}$. Then from
(\ref{eq3}) we get: $\mathscr{W}_m\subseteq \mathscr{Y}_{m+d}$;
from (\ref{eq2}) we get: if $J\in\mathscr{W}_m$ and
$K\in\mathscr{Y}_m$, then $J\cap K\in\mathscr{Y}_{m+d}$; and from
(\ref{eq4}) we get: if $J,K\in\mathscr{Y}_{m}$ and $J\cap
K=\emptyset$, then $J\sqcup K\in\mathscr{Y}_{m+1}$.

By induction, we deduce $\mathcal{I}_k(\mathscr{W}_m)\subseteq
\mathscr{Y}_{m+kd}$ for all $k$, and
$\mathcal{V}_k(\mathscr{Y}_m)\subseteq\mathscr{Y}_{m+k}$ for all
$k$. Composing, we obtain
$\mathcal{V}_k(\mathcal{I}_l(\mathscr{W}_m))\subseteq
\mathcal{V}_k(\mathscr{Y}_{m+ld})\subseteq\mathscr{Y}_{m+ld+k}$.
By 2) of Lemma \ref{RnD2n}, we get
$\mathcal{D}^n(\mathscr{W}_m)\subseteq\mathscr{Y}_{m+2^nd+2^{2^n}}$.
Hence, using (\ref{Eq5}), we obtain
$\mathscr{P}(I)=\mathscr{Y}_D$, where $D=M+4^Nd+2^{4^N}$.

Let $B$ denote the subgroup generated by $b$, so that $N$ is the
normal subgroup generated by $B$. Let $r$ be the order of $b$.
Then $B^I$ is contained in $X_{D+r}$. Moreover, there exists $R$
such that every element of $N$ is the product of $R$ conjugates of
elements of $B$. Then, using that $c_I\in X_0$ for all $c\in G$,
$N^I$ is contained in $X_{D+r+3R}$.~$\blacksquare$

\begin{thm}
Let $G$ be a finite group, and let $N$ the last term of its
descending central series (so that $[G,N]=N$). Let $I$ be any set,
and set $H=G^I$. Let $(X_m)$ be an increasing sequence of subsets
of $H$ such that $\mathcal{G}(X_m)\subseteq X_{m+1}$ and $\bigcup
X_m=H$. Then $N^I\subseteq X_m$ for $m$ big
enough.\label{Last_term}
\end{thm}

\noindent \textbf{Proof}:~Let $G$ be a counterexample with $|G|$
minimal. Let $W$ be a normal subgroup of $G$ such that $W^I$ is
contained in $X_m$ for large $m$, and which is maximal for this
property. Since $G$ is a counterexample, $N\nsubseteq W$. Hence
$G/W$ is not nilpotent, and is another counterexample, so that, by
minimality, $W=\{1\}$. Since $G$ is not nilpotent, there exists,
by Lemma \ref{f(a,b)}, $a\in G$, $b\in G-\{1\}$, and a homogeneous
monomial $f:G^2\to G$, such that $f(a,b)=b$. So, if $M$ is the
normal subgroup generated by $b$, $M^I$ is contained, by Lemma
\ref{Main_Lemma}, in $X_i$ for large $i$. This contradicts the
maximality of $W\;(=\{1\})$.~$\blacksquare$

\medskip

In view of Proposition \ref{strong_bound}, Theorem
\ref{G^I_strongly_bounded} immediately follows from Theorem
\ref{Last_term}. Theorem \ref{G^I_strongly_bounded} has been
independently proved by Khelif \cite{Khelif}, who also proves
Proposition \ref{cof_P(I)}, but concludes by another method.

\begin{que}
Let $G$ be a finite group, and $N$ a subgroup of $G$ which
satisfies the conclusion of Theorem \ref{Last_term} ($I$ being
infinite). Is it true that, conversely, $N$ must be contained in
the last term of the descending central series of $G$? We
conjecture that the answer is positive, but the only thing we know
is that $N$ must be contained in the derived subgroup of
$G$.\label{quest_nilp}
\end{que}

\begin{rem}
We could have introduced a relative definition: if $G$ is a group
and $X\subseteq G$ is a subset, we say that the pair $(G,X)$ is
strongly bounded if, for every isometric action of $G$ on any
metric space $M$ and every $m\in M$, then the ``$X$-orbit" $Xm$ is
bounded. Note that $G$ is strongly bounded if only if the pair
$(G,G)$ is strongly bounded. Proposition \ref{strong_bound}
generalizes as: the pair $(G,X)$ is strongly bounded if and only
if for every sequence $(X_n)$ of subsets of $G$ such that
$\bigcup_n X_n=G$ and $\mathcal{G}(X_{n})\subseteq X_{n+1}$ for
all $n$, one has $X_n\supseteq X$ for some~$n$.

Theorem \ref{Last_term} is actually stronger than Theorem
\ref{G^I_strongly_bounded}: it states that if $G$ is a finite
group, if $N$ is the last term of its descending central series,
and if $I$ is any set, then the pair $(G^I,N^I)$ is strongly
bounded. This provides non-trivial strongly bounded pairs of
solvable groups (trivial pairs are those pairs $(G,X)$ with $X$
finite); compare Remark \ref{rem:solv_SB} and Question
\ref{quest_nilp}.
\end{rem}


\begin{que}
We do not assume the continuum hypothesis. Does there exist a
strongly bounded group with cardinality $\aleph_1$?
\end{que}
It seems likely that a variation of the argument in \cite{Shelah}
might provide examples.

\begin{que}Let $(G_n)$ be a sequence of finite perfect groups. When
is the product $\prod_{n\in\mathbf{N}}G_n$ strongly bounded?
\end{que}
It follows from Theorem \ref{G^I_strongly_bounded} that if the
groups $G_n$ have bounded order, then $\prod_{n\in\mathbf{N}}G_n$
is strongly bounded. If all $G_n$ are simple, Saxl, Shelah and
Thomas prove \cite[Theorems 1.7 and 1.9]{SST} that
$\prod_{n\in\mathbf{N}}G_n$ has cofinality $\neq\omega$ if and
only if there does \textit{not} exist a fixed (possibly twisted)
Lie type $L$, a sequence $(n_i)$ and a sequence $(q_i)$ of prime
powers tending to infinity, such that $G_{n_i}\simeq L(q_i)$ for
all~$i$. Does this still characterize infinite strongly bounded
products of non-abelian finite simple groups?

\appendix

\section{Proof of Lemma \ref{f(a,b)}}

This Appendix is added for the convenience of the reader. It is
dropped in the published version.

\begin{lem}[\cite{KT}]
Let $G$ be a group, $g\in G$, and $g'$ an element of the subgroup
generated by the conjugates of $g$. Then there exists a
homogeneous monomial $f:G\to G$ such that
$f(g)=g'$.\label{exist_hom_mon}
\end{lem}
\noindent \textbf{Proof}:~Write $g'=\prod
c_ig^{\alpha_i}c_i^{-1}$. Then $x\mapsto \prod
c_ix^{\alpha_i}c_i^{-1}$ is a homogeneous monomial and
$f(g)=g'$.~$\blacksquare$

\medskip


\begin{lem}
Let $G$ be a finitely generated group. Suppose that $G$ is not
nilpotent. Then there exists $a\in G$ such that the normal
subgroup of $G$ generated by $a$ is not nilpotent.\label{nnilp}
\end{lem}
\noindent \textbf{Proof}:~Fix a finite generating subset $S$ of
$G$. For every $s\in S$, denote by $N_s$ the normal subgroup of
$G$ generated by~$s$. Since finitely many nilpotent normal
subgroups generate a nilpotent subgroup, it immediately follows
that if all $N_s$ are nilpotent, then $G$ is
nilpotent.~$\blacksquare$

\medskip

\noindent \textbf{Proof} of Lemma \ref{f(a,b)}. We reproduce the
proof from \cite{KT}. Let $G$ be a finite group which is not
nilpotent. We must show that there exist $a\in G$, $b\in G-\{1\}$,
and a homogeneous monomial $f:G^2\to G$, such that $f(a,b)=b$.

Take $a$ as in Lemma \ref{nnilp}, and $A$ the normal subgroup
generated by $a$. Let $A_1$ be the upper term of the ascending
central series of $A$. We define inductively the sequences
$(a_i)_{i\in\mathbf{N}}$ and $(b_i)_{i\in\mathbf{N}}$ such that
$$b_i\in A-A_1,\qquad a_i\in A\qquad\text{and}\qquad
b_{i+1}=[a_i,b_i]\in A-A_1.$$ Since $G$ is finite, there exist
integers $m,m'$ such that $m<m'$ and $b_m=b_{m'}$. Set $b=b_m$,
and for all $i$, choose, using Lemma \ref{exist_hom_mon}, a
homogeneous monomial $f_i$ such that $f_i(a)=a_i$. Then the
monomial
$$f:\;\;\;(x,y)\mapsto
[f_{m'-1}(x),[f_{m'-2}(x),\dots,[f_m(x),y],\dots]]$$ satisfies
$f(a,b)=b$.~$\blacksquare$

\section{Groups with cardinality $\aleph_1$ and Property (FH)}

This appendix is dropped in the published version.

\begin{prop}
Let $G$ be a countable group. Then $G$ embeds in a group of
cardinality $\aleph_1$ with Property (FH).\label{aleph_1_FH}
\end{prop}

The proof rests on two ingredients.

\begin{thm}[Delzant]
If $G$ is any countable group, then $G$ can be embedded in a group
with Property (T).\label{Delzant}
\end{thm}
{\bf Sketch of proof}: this is a corollary of the following
result, first claimed by Gromov\footnote{Theorem 5.6.E in {\em
Hyperbolic groups}. In ``Essays in group theory" (S. Gersten,
ed.), MSRI series vol. 8, Springer Verlag, 1987.}, and
subsequently independently proved by Delzant\footnote{{\em
Sous-groupes distingu\'es et quotients des groupes hyperboliques}.
\newblock Duke Math. J., {\bf 83}, Vol. 3, 661-682, 1996.} and
Olshanskii\footnote{{\em SQ-universality of hyperbolic groups},
Sbornik Math. {\bf 186}, no. 8, 1199-1211, 1995.}: if $H$ is any
non-elementary word hyperbolic group, then $H$ is SQ-universal,
that is, every countable group embeds in a quotient of $H$. Thus,
the result follows from the stability of Property (T) by
quotients, and the existence of non-elementary word hyperbolic
groups with Property (T); for instance, uniform lattices in
$\textnormal{Sp}(n,1)$, $n\ge 2$ (see \cite{HV}).~$\blacksquare$



\medskip

Let $\mathcal{C}$ be any class of metric spaces, let $G$ be a
group. Say that $G$ has Property (F$\mathcal{C}$) if for every
isometric action of $G$ on a space $X\in\mathcal{C}$, all orbits
are bounded. For instance, if $\mathcal{C}$ is the class of all
Hilbert spaces, then we get Property (FH).

\begin{prop}
Let $G$ be a group in which every countable subset is contained in
a subgroup with Property (F$\mathcal{C}$). Then $G$ has Property
(F$\mathcal{C}$).
\label{indFH}
\end{prop}
\noindent \textbf{Proof}:~Let us take an affine isometric action
of $G$ on a metric space $X$ in $\mathcal{C}$, and let us show
that it has bounded orbits. Otherwise, there exists $x\in X$, and
a sequence $(g_n)$ in $G$ such that $d(g_n\,x,x)\to\infty$. Let
$H$ be a subgroup of $G$ with Property (F$\mathcal{C}$) containing
all $g_n$. Since $Hx$ is not bounded, we have a
contradiction.~$\blacksquare$

\medskip

\noindent \textbf{Proof} of Proposition \ref{aleph_1_FH}. We make
a standard transfinite induction on $\omega_1$ (as in
\cite{Sabbagh}), using Theorem \ref{Delzant}. For every countable
group $\Gamma$, choose a proper embedding of $\Gamma$ into a group
$F(\Gamma)$ with Property (T) (necessarily finitely generated).
Fix $G_0=G$, $G_{\alpha+1}=F(G_\alpha)$ for every
$\alpha<\omega_1$, and
$G_\lambda=\underrightarrow{\lim}_{\beta<\lambda}G_\lambda$ for
every limit ordinal $\lambda\le\omega_1$. It follows from
Proposition \ref{indFH} that $G_{\omega_1}$ has Property (FH).
Since all embeddings $G_\alpha\to G_{\alpha+1}$ are proper,
$G_{\omega_1}$ is not countable, hence has cardinality
$\aleph_1$.~$\blacksquare$

\bigskip

\noindent\textbf{Acknowledgments.} I thank Bachir Bekka, who has
suggested to me to show that the groups studied in \cite{KT} have
Property (FH). I am grateful to George M. Bergman, David Madore,
Romain Tessera and Alain Valette for their useful corrections and
comments.

\bigskip
\footnotesize
\noindent Yves de Cornulier\\
E-mail: \url{decornul@clipper.ens.fr}\\
\'Ecole Polytechnique F\'ed\'erale de Lausanne (EPFL)\\
Institut de G\'eom\'etrie, Alg\`ebre et Topologie (IGAT)\\
CH-1015 Lausanne, Switzerland


\begin{thebibliography}{KM98b}

\bibitem[BHV]{BHV} Bachir {\sc Bekka}, Pierre {\sc de la Harpe}, Alain {\sc Valette}.
\newblock ``Kazhdan's Property (T)". \newblock Forthcoming book, currently available at
\url{http://poncelet.sciences.univ-metz.fr/~bekka/}, 2004.

\bibitem[BH]{BH} Martin R. {\sc Bridson}, Andr\'e {\sc Haefliger}.
\newblock ``Metric Spaces of Non-Positive Curvature". \newblock
Grundlehren Math. Wiss. 319, Springer, 1999.

\bibitem[Ber04]{Bergman} George M. {\sc Bergman}.
\newblock {\em Generating infinite symmetric groups}, Bull. London
Math. Soc., to appear; arXiv math.GR/0401304, 2005.


\bibitem[DH05]{DH} Manfred {\sc Droste}, W. Charles {\sc Holland}.
\newblock {\em Generating automorphism groups of chains}.
\newblock Forum Math. {\bf 17}, 699-710, 2005.

\bibitem[DG05]{DG} Manfred {\sc Droste}, R\"udiger {\sc G\"obel}.
\newblock {\em Uncountable cofinalities of permutations groups}.
\newblock J. London Math. Soc. {\bf 71}(2), 335-344, 2005.

\bibitem[HV]{HV} Pierre {\sc de la Harpe}, Alain {\sc Valette}.
\newblock ``La propri\'et\'e (T) de Kazhdan pour les groupes localement
compacts", Ast\'erisque {\bf 175}, SMF, 1989.

\bibitem[Khe05]{Khelif} Anatole {\sc Khelif}. \newblock {\em \`A
propos de la propri\'et\'e de Bergman}. \newblock Preprint 2005.

\bibitem[LS]{LS} Roger C. {\sc Lyndon}, Paul E. {\sc Schupp}. \newblock ``Combinatorial
group theory", Springer, 1977.

\bibitem[KR05]{KR} Alexander S. {\sc Kechris}, Christian {\sc
Rosendal}. \newblock {\em Turbulence, amalgamation and generic
automorphisms of homogeneous structures}. \newblock Preprint,
2005. ArXiv math.LO/0409567.

\bibitem[KT74]{KT} Sabine {\sc Koppelberg}, Jacques {\sc Tits}.
\newblock {\em Une propri\'et\'e des produits directs infinis de
groupes finis isomorphes}. \newblock C. R. Math. Acad. Sci. Paris,
S\'er. A {\bf 279}, 583-585, 1974.

\bibitem[Mil04]{Miller} Benjamin D. {\sc Miller}. \newblock {\em Full groups,
classification, and equivalence relations}. \newblock PhD
dissertation, 2004. Currently available at
\url{http://www.math.ucla.edu/~bdm/papers.html}.

\bibitem[Sab75]{Sabbagh} Gabriel {\sc Sabbagh}. \newblock {\em Sur les groupes qui ne
sont pas r\'eunion d'une suite croissante de sous-groupes
propres}.
\newblock C. R. Math. Acad. Sci. Paris, S\'er. A {\bf 280}, 763-766, 1975.

\bibitem[Sco51]{Scott} William R. {\sc Scott}. \newblock {\em Algebraically closed
groups}. \newblock Proc. Amer. Math. Soc. {\bf 2}, 118-121, 1951.

\bibitem[Ser]{Serre} Jean-Pierre {\sc Serre}. \newblock
``Arbres, amalgames, $\text{SL}_2$". \newblock Ast\'erisque {\bf
46}, SMF, 1977.

\bibitem[She80]{Shelah} Saharon {\sc Shelah}. \newblock {\em On a problem of Kurosh,
J\'onsson groups, and applications}. \newblock In {\em Word
problems II}. \newblock North-Holland Publ. Company, 373-394,
1980.

\bibitem[SST96]{SST} Jan {\sc Saxl}, Saharon {\sc Shelah}, Simon
{\sc Thomas}. \newblock {\em Infinite products of finite simple
groups}. \newblock Trans. Amer. Math. Soc. {\bf 348}(11),
4611-4641, 1996.


\bibitem[W01]{W} Report of the workshop {\em Geometrization of
Kazhdan's Property (T)} (organizers: B. Bekka, P. de la Harpe, A.
Valette; 2001).
\newblock
Unpublished; currently available at
\url{http://www.mfo.de/cgi-bin/tagungsdb?type=21&tnr=0128a}.

\end{thebibliography}
\end{document}